\documentclass[12pt]{amsart}

\usepackage{euscript,amscd,enumerate,hyperref}

\usepackage{comment}

\usepackage{mathrsfs}

\usepackage{graphicx}

\usepackage{color}

\usepackage[margin=0.8in]{geometry}

\usepackage{bbm}

\newtheorem{theorem}{Theorem}[section]
\newtheorem*{theorem*}{Theorem}

\newtheorem{assumption}[theorem]{Assumption}

\theoremstyle{remark}
\newtheorem{rmk}[theorem]{Remark}

\newcommand{\E}{\mathbb E}

\newcommand{\rr}{\mathbb{R}}

\newcommand{\nn}{\mathbb{N}}

\newcommand{\eq}{\begin{equation}}
\newcommand{\en}{\end{equation}}

\numberwithin{equation}{section}

\title{Large rank-based models with common noise}

\author{Praveen Kolli}
\address{Department of Mathematical Sciences, Carnegie Mellon University}
\email{kpc.0915@gmail.com}

\author{Andrey Sarantsev}
\address{Department of Mathematics and Statistics, University of Nevada, Reno}
\email{asarantsev@unr.edu}

\begin{document}

\thanks{The first author thanks his adviser Mykhaylo Shkolnikov for many stimulating conversations, and for pointing out several mistakes in the earlier versions of the paper. The first author thanks Julien Reygner for many helpful conversations. Finally, the first author thanks Ioannis Karatzas and Johannes Ruf for posing problems similar to the one discussed in this paper.}



\begin{abstract}
For large systems of Brownian particles interacting through their ranks introduced in (Banner, Fernholz, Karatzas, 2005), the empirical cumulative distribution function satisfies a porous medium PDE. However, when we introduce a common noise, the limit is no longer deterministic. Instead, we show that this limit is a solution of a stochastic PDE related to this porous medium PDE. This stochastic PDE is somewhat similar to the equations developed for conservation laws with rough stochastic fluxes (Lions, Perthame, Souganidis, 2013). 
\end{abstract}

\maketitle

\thispagestyle{empty}

\section{Introduction}

We study interacting particle systems $X_1^{(n)}(t), \ldots, X_n^{(n)}(t),\, t \ge 0$,  on the real line,    governed by the following system of stochastic differential equations: 
\begin{equation}\label{main_sde}
\begin{split}
\mathrm{d}X^{(n)}_i(t) &= b\big(F_{\rho^{(n)}(t)}\big(X^{(n)}_i(t)\big)\big)\,\mathrm{d}t + \sigma\big(F_{\rho^{(n)}(t)}\big(X^{(n)}_i(t)\big)\big)\,\mathrm{d}B^{(n)}_i(t) \\ & \quad +\,\gamma\big(t,\rho^{(n)}(t)\big)\mathrm{d}W(t),\, i = 1, \ldots, n;
\\ \rho^{(n)}(t) & :=\frac{1}{n}\sum_{i=1}^n \delta_{X^{(n)}_i(t)}.
\end{split}
\end{equation}
Here, $F_{\mu}$ is the cumulative distribution function of a probability measure $\mu$, and $b : [0, 1] \to (0, \infty)$, $\sigma : [0, 1] \to (0, \infty)$ are given functions. We fix a $p > 1$ and denote by $\mathcal P_p$ the space of all probability measures on 
$\mathbb R$ with finite $p$th moment, so $\gamma :  [0, \infty) \times \mathcal P_p(\rr) \to \rr$ is another given function. Finally, $W,\,B^{(n)}_1,\,B^{(n)}_2,\,\ldots,\,B^{(n)}_n$, are i.i.d. standard Brownian motions. 

In the absence of the common noise: $\gamma \equiv 0$, the system \eqref{main_sde} reduces to a rank-based model. Originally, rank-based model appeared as a special case in the context of the piecewise linear filtering problem in \cite{BP} where weak uniqueness is established (weak existence being a consequence of the general result in \cite[Exercise 12.4.3]{SV}). Rank-based models have attracted a lot of attention recently since their appearance in stochastic portfolio theory where they are used to model stock prices in large equity markets (\cite[section 13]{FK}, \cite{BFK}, \cite{IPS}, \cite{CP}). Let $S_i(t)$ denote the market capitalization (the number of shares, multiplied by the share price) of the $i$th company, $i = 1, \ldots, n$, listed in any of the major  stock exchanges. Then $(\log S_1(t), \ldots, \log S_n(t))$ is modeled as a rank-based model.  

A limitation of these rank-based models used for large equity markets is that independent Brownian motions drive stock prices (even though rank-based particles themselves do not evolve independently). A richer model, with correlated Brownian motions driving the rank-based particles, would probably better capture the characteristics of a large equity market. One could consider various covariance matrices, for example block-diagonal matrices, corresponding to correlated fluctuations for stocks in the same economy sector, and stocks from different sectors having independent Brownian fluctuations. In fact, estimating the covariance matrix (and its inverse, which is called the precision matrix) in high dimension is a well-developed research field, see for example the book \cite{Cov} and the survey \cite{Fan}. To apply this theory for rank-based stock market models is left for future research. 

The model~\eqref{main_sde} is a first step in this direction: It has only one Brownian term common to all stocks, and another independent Brownian term for each stock. 

In this paper, we are concerned with the large system limit (as the number $n$ of particles tends to infinity) for the particle system in \ref{main_sde}. This can be viewed as a version of a law of large numbers. Such a result would give us an understanding of the behavior of the whole market under the assumption that the number of firms operating in the market is large. 

We state some known results about the limiting behavior for large systems from~\eqref{main_sde} with $\gamma \equiv 0$ (that is, classic rank-based systems, without common noise). Under suitable regularity conditions, it was shown  in \cite[Proposition 2.1]{JR} that the measure-valued processes $\mu^{(n)}(\cdot)$, $n\in\nn$ (see \eqref{eq:Y} below and the discussion following it) converge (in a certain functional space, see below) to a deterministic limit $\mu(\cdot)$, whose cumulative distribution function $R(t,\cdot):=F_{\mu(t)}(\cdot)$ satisfies the \textit{porous medium equation}:   
\begin{equation}\label{PME}
R_t=-B(R)_x+\Sigma(R)_{xx},\quad R(0,\cdot)=F_\lambda(\cdot),
\end{equation}
in the generalized sense, see \cite[Definition 3]{Gi}, with 
$$
B(r):=\int_0^r b(a)\,\mathrm{d}a\quad \mbox{and}\quad \Sigma(r):=\int_0^r \frac{1}{2}\,\sigma(a)^2\,\mathrm{d}a.
$$
A central limit theorem type result was obtained in \cite[Theorem 1.2]{KoS}: fluctuations of the empirical cumulative distribution function around its limit are governed by a suitable SPDE; and a large deviations result was obtained in \cite[Theorem 1.4]{DSVZ}.

Since there is common noise in~\eqref{main_sde}, we cannot hope that as $n \to \infty$, all noise will be canceled. The limit $\rho$, in other words, will be stochastic: it satisfies the SPDE~\eqref{SPDE2} below.  Let $\mathcal H^{l}\big(\rr)$ denote the H\"older space, see \cite[Pg. 7]{LSU}. For a $p \ge 1$ and a metric space $E$, $\mathcal P_p(E)$ is the space of probability measures on $E$ with finite $p$th moment,  equipped with {\it Wasserstein distance} $\mathcal W_p$, which is defined for $\mu, \nu \in \mathcal P_p(\rr)$ as follows:
$$
\mathcal W_p(\mu,\nu)=\inf_{(Y_1,Y_2)} \E\big[|Y_1-Y_2|^p\big]^{1/p},
$$
with the $\inf$ over random vectors $(Y_1,Y_2)$ such that $Y_1 \sim \mu$ and $Y_2 \sim \nu$.  Also, $C([0, T], E)$ is the space of continuous functions $[0, T] \to E$ with sup-distance, and $D([0, T], E)$ stands for the Skorohod space of right-continuous functions with left limits from $[0, T]$ to $E$. 

\smallskip

We fix a $p > 1$ and a time horizon $T > 0$.

\begin{assumption}\label{main_asmp}
\begin{enumerate}[(a)]
\item The functions $b,\, \sigma$  are differentiable, and $b', \sigma' \in \mathcal H^{\beta}(\mathbb R)$ for a $\beta > 0$.
\item The function $\sigma$ is bounded away from zero: $\min_{a\in[0,1]} \sigma(a)>0$.
\item The function $\gamma$ is bounded and Lipschitz with respect to $\mathcal W_1$:
\begin{equation}
\label{eq:gamma-lip}
|\gamma(t, \nu_1) - \gamma(t, \nu_2)| \le L_{\gamma}\mathcal W_1(\nu_1, \nu_2),\ \mbox{for all}\ \ t \in [0, T]\ \mbox{and}\ \nu_1, \nu_2 \in \mathcal P_1.
\end{equation}
\item There exists a measure $\lambda \in \mathcal P_p$ such that $X^{(n)}_1(0),X^{(n)}_2(0),\ldots,X^{(n)}_n(0) \sim \lambda$ i.i.d
\item The cumulative distribution function of $\lambda$ satisfies $F_\lambda(\cdot) \in \mathcal H^{3 + \beta}\big(\rr)$.
\end{enumerate}
\end{assumption}

We are now ready to state the main result of the paper.

\begin{theorem}\label{thm:lln}
Under Assumption~\ref{main_asmp}, for every $n$ the system~\eqref{main_sde} has a unique solution in the weak sense, unique in law. For every $T > 0$ and $q \in [1, p)$, the sequence $(\rho^{(n)})_{n \ge 1}$ of random elements in   $C([0,T], \mathcal P_q)$ weakly converges to $\rho$, a unique solution in $C([0, T], \mathcal P_1)$ to the following functional equation:
\begin{equation}\label{main_eqn} 
F_{\rho(t)}(x)=R(t,x-\Gamma(t)),\quad \Gamma(t) := \int_0^t \gamma(\rho(s))\mathrm{d}W(s).
\end{equation}
The function $G(t, \cdot) := F_{\rho(t)}(\cdot)$ solves the following SPDE: 
\begin{equation}\label{SPDE2}
\mathrm{d}G=\big[-B(G)_x+\Sigma(G)_{xx}+\frac{1}{2}G_{xx}\gamma^2(t, \rho(t))\big]\mathrm{d}t\,-\gamma\big(t, \rho(t)\big)G_x \mathrm{d}W(t).
\end{equation}
\end{theorem}

\begin{rmk}
The SPDE in \eqref{SPDE2} is very closely related to stochastic scalar conservation laws introduced by Lions, Perthame and Souganidis, \cite{LPS}. They introduced the notion of pathwise entropy solutions to stochastic scalar conservation laws and this  theory was extended to a certain class of problems in \cite{GS}. In particular, if $b=0$ and $\gamma=1$, the SPDE in \eqref{SPDE2} reduces to the SPDE in \cite[Equation 1.1]{GS} with $F(x)=x$ and $A(u)=\frac{1}{2}\,\sigma(u)^2$ and in this case $G(t,x) = R(t,x-W(t))$ is a solution of the SPDE. This observation opens the door to further research on stochastic scalar conversation laws from the perspective of rank-based models with common noise.
\end{rmk}

\begin{rmk}
The functional equation in \eqref{main_eqn} admits an explicit representation when the function $\gamma(t,\nu)=f(t)$, where $f$ is any continuous function defined on $[0,\infty)$. Another case of special interest is the function $\gamma(t,\nu)=\int_{\rr}f(x)\nu(\mathrm{d}x)$, where the function $f$ defined on $\rr$ is differentiable and has a bounded derivative. Integrating by parts, we get:
\begin{equation}
\int_{\rr}f(x)\nu(\mathrm{d}x)=-\int_{\rr}f'(x)F_{\nu}(x)\mathrm{d}x.
\end{equation}
Using representation of $\mathcal W_1$ from \cite[p.64]{SW}, we prove the Lipschitz property of $\gamma$:
\begin{align*}
|\gamma(t, \nu_1) & - \gamma(t, \nu_2)| = \Bigl|\int_{\rr}f'(x)\left[F_{\nu_1}(x) - F_{\nu_2}(x)\right]\,\mathrm{d} x\Bigr| \\ & \le 
\sup_{x \in \rr}|f'(x)|\int_{\rr}\left[F_{\nu_1}(x) - F_{\nu_2}(x)\right|\,\mathrm{d} x = \sup_{x \in \rr}|f'(x)|\cdot \mathcal W_1(\nu_1, \nu_2). 
\end{align*}
\end{rmk}

\section{Proof of Theorem~\ref{thm:lln}}

\subsection{Overview of the proof} For notational convenience, we assume $\gamma$ does not depend on the time variable. 
We split the proof into four sections. In the first, we simply prove existence and uniqueness for the finite system~\eqref{main_sde}. In the second section, we prove tightness of the sequence $(\rho^{(n)})$ of empirical measures in $D([0, T], \mathcal P_q)$. In the third section, we prove that every weak limit point $\rho$ has the cumulative distribution function $F_{\rho(t)}$ solving~\eqref{main_eqn}. In the fourth section, we derive~\eqref{SPDE2} from~\eqref{main_eqn}. Finally, in the fifth section we prove uniqueness for solutions of~\eqref{main_eqn}.

\subsection{Existence and uniqueness of the finite system} First,  let us show weak existence and uniqueness in law of the system~\eqref{main_sde}. Define
\begin{equation}
\label{eq:Y}
Y^{n}_i(t)=X^{(n)}_i(t)-\int_0^t \gamma\big(\rho^{(n)}(s)\big)\mathrm{d}W(s), \quad \mu^{(n)}(t):=\frac{1}{n}\sum_{i=1}^n \delta_{Y^{(n)}_i(t)}.
\end{equation}
It is straightforward to check that 
\begin{equation}\label{observations}
\begin{split}
F_{\rho^{(n)}(t)}(x) & = F_{\mu^{(n)}(t)}\Big(x-\int_0^t \gamma\big(\rho^{(n)}(s)\big)\mathrm{d}W(s)\Big);\\ F_{\rho^{(n)}(t)}\big(X^{(n)}_i(t)\big) & = F_{\mu^{(n)}(t)}\big(Y^{(n)}_i(t)\big).\\
\end{split}
\end{equation}
Therefore,~\eqref{eq:Y} satisfy the system of equations similar to~\eqref{main_sde}, but with $\gamma = 0$. This is our key observation. The classic system of competing Brownian particles $Y^{(n)} = (Y^{(n)}_1, \ldots, Y^{(n)}_n)$ exists in the weak sense and is unique in law. We can rewrite~\eqref{eq:Y} as follows:
\begin{equation}
\label{eq:fixed-point}
X^{(n)}_i(t) = Y^{(n)}_i(t) +  \int_0^t\gamma\big(\rho^{(n)}(s)\big)\mathrm{d}W(s),\quad 0 \le t \le T. 
\end{equation}
Define the space $\mathfrak X$ of all random elements on our filtered probability space with values in the metric space $C([0, T], \mathcal P_1(\mathbb R))$, and with finite second moment. 
Define the mapping $\Phi: \mathfrak X \to \mathfrak X$ as follows: Fix $\mu \in \mathcal X$ and fix some realization $\mu(\omega)$. For $t \in [0, T]$, $\Phi(\mu)(t)$ is the empirical distribution of $n$ particles 
$$
X^{(n, \mu)}_i(t) := Y^{(n)}_i(t) +  \int_0^t\gamma\big(\mu(s)\big)\mathrm{d}W(s),\, i = 1, \ldots, n,
$$
Then~\eqref{eq:fixed-point} is equivalent to saying that $\rho^{(n)}$ is a fixed point of the mapping $\Phi$. Couple $\Phi(\mu)(t)$ and $\Phi(\nu)(t)$ as the uniform measure on the set $\{(X_i^{(n, \mu)}(t), X_i^{(n, \nu)}(t))\mid i = 1, \ldots, n\}$ (if some of these points coincide, we count them twice in this measure). From properties of the It\^o integral and Lipschitz condition~\eqref{eq:gamma-lip}, 
\begin{align}
\label{eq:major-estimate}
\begin{split}
\mathbb E\,&\mathcal W^2_1(\Phi(\mu)(t), \Phi(\nu)(t))  \le \frac1n\sum\limits_{i=1}^n\mathbb E\,\bigl|X^{(n, \mu)}_i(t) - X^{(n, \nu)}_i(t)\bigr|^2 \\ & \le \frac1n\sum\limits_{i=1}^n\mathbb E\,\Bigl|\int_0^t\{\gamma(\mu(u)) - \gamma(\nu(u))\}\,\mathrm{d}W(u)\Bigr|^2 \\ & = \int_0^t\,\mathbb E\,(\gamma(\mu(u)) - \gamma(\nu(u)))^2\,\mathrm{d}u \le L_{\gamma}^2\int_0^t\,\mathbb E\,\mathcal W^2_1(\mu(s), \nu(s))\,\mathrm{d}s.
\end{split}
\end{align}
Iterating $\Phi$ and integrating by parts, we get similarly to 
\cite[Section 5.2.B, (2.19)]{KS}
$$
\mathbb E\sup\limits_{0 \le s \le t}\mathcal W^2_1(\Phi^k(\mu)(s), \Phi^k(\nu)(s)) \le 
\frac{(4L_{\gamma}^2)^k}{k!}\mathbb E\sup\limits_{0 \le s \le t}\mathcal W^2_1(\mu(s), \nu(s)).
$$
The rest of the proof is the standard argument, see \cite[Section 5.2.B, (2.19)]{KS}: We fix $\nu_0 \in \mathcal P_1$, use the Borel-Cantelli lemma to prove convergence of the sequence $(\Phi^n\nu_0)$ a.s. in $C([0, T], \mathcal P_1)$ to some $\nu$. This limit satisfies $\Phi\nu = \nu$. From~\eqref{eq:major-estimate}, we get uniqueness of the fixed point for $\Phi$. This completes the proof of existence and uniqueness.


\subsection{Tightness} Let us show that the sequence $(\rho^{(n)})_{n \ge 1}$ is tight in $D([0, T], \mathcal P_q)$ for $q < p$. We follow the proof of \cite[Lemma 7.4]{MyOwn19}. Apply It\^o's formula to $\bigl(\rho^{(n)}_t, f\bigr)$ for $f \in C^2_b(\mathbb R)$ (the space of $C^{2}$ functions $\mathbb R \to \mathbb R$  bounded together with their first and second derivatives): 
\begin{align}
\label{eq:big-ito}
\begin{split}
\mathrm{d}& \bigl(\rho^{(n)}_t, f\bigr)  = \frac1n\sum\limits_{k=1}^nf'\bigl(X_i^{(n)}(t)\bigr)b\bigl(X_i^{(n)}(t)\bigr)\,\mathrm{d}t \\ & + \frac1n\sum\limits_{k=1}^nf'\bigl(X_i^{(n)}(t)\bigr)\sigma\bigl(X_i^{(n)}(t)\bigr)\,\mathrm{d} B_i(t) + \frac1n\sum\limits_{k=1}^nf'\bigl(X_i^{(n)}(t)\bigr)\gamma\bigl(\rho^{(n)}_t\bigr)\,\mathrm{d}W(t) \\ & + 
\frac{1}{2n}\sum\limits_{i=1}^nf''\bigl(X_i^{(n)}(t)\bigr)\sigma^2\bigl(X_i^{(n)}(t)\bigr)\,\mathrm{d} t + \frac{1}{2n}\sum\limits_{i=1}^nf''\bigl(X_i^{(n)}(t)\bigr)\gamma^2\bigl(\rho^{(n)}_t\bigr)\,\mathrm{d}t.
\end{split}
\end{align}
Since $f', f'', b, \sigma, \gamma$ are bounded, the equation~\eqref{eq:big-ito} has the form 
$$
\mathrm{d} \bigl(\rho^{(n)}_t, f\bigr) = \alpha_n(t)\,\mathrm{d} t + \theta_n(t)\,\mathrm{d}\tilde{B}_n(t),\, n = 1, 2, \ldots;\, t \in [0, T],
$$
for uniformly bounded $\alpha_n$ and $\theta_n$, and for Brownian motions $\tilde{B}_n$. Thus application of standard tools gives us tightness of $(\rho^{(n)}_t, f)_{n \ge 1}$ in $C[0, T]$, and therefore in $D[0, T]$. From the Burkholder-Davis-Gundy inequality  \cite[Theorem 3.28]{KS} and boundedness of $b$, $\sigma, \gamma$, we get: 
$$
\mathbb E\Bigl[\max\limits_{0 \le t \le T}\bigl|X_i^{(n)}(t)\bigr|^p\Bigr] \le  C < \infty.
$$
Thus, the sequence of measure-valued processes $(\rho^{(n)})$ satisfies for $f_p(x) := |x|^p$:
\begin{equation}
\label{eq:p-bound}
\mathbb E \sup\limits_{0 \le t \le T}\bigl(\rho^{(n)}_t, f_p\bigr) \le C.
\end{equation}
Take any $\eta > 0$, and consider the subset $\mathcal K := \{\nu \in \mathcal P_q\mid (\nu, f_p) \le C/\eta\}$, which is compact in $\mathcal P_{q}$ by \cite[Lemma 2.2]{MyOwn19}. From the standard Markov inequality, we have:
$$
\mathbb P\Bigl[\rho_t^{(n)} \in \mathcal K\quad \forall\, t \in [0, T]\Bigr] > 1 - \eta.
$$
Next, take the algebra $\mathfrak A$ in $C_b(\mathcal P_q)$, the space of bounded continuous functions $\mathcal P_q \to \mathbb R$, generated by $\mathfrak M := \{(\cdot, f)\mid f \in C^2_b\}$. This set $\mathfrak M$ separates points: for every $\nu'$ and $\nu''$ in $\mathcal P_q$, there exists an $f \in C^2_b$ such that $(\nu', f) \ne (\nu'', f)$. This set $\mathfrak M$ also contains $1$, because $f_0 = 1 \in C^2_b$. By the Stone-Weierstrass theorem \cite[Section 4.7]{Folland}, the algebra $\mathfrak A$ is dense in $C_b(\mathcal P_q)$ in the topology of uniform convergence on compact subsets. Note that $(\rho^{(n)}_t, f)$ is uniformly bounded for $f \in C_b^2(\mathbb R)$. Therefore, for every collection $g_1, \ldots, g_m \in C_b^2(\mathbb R)$, the following sequence is tight in $C[0, T]$ (and therefore in $D[0, T]$, \cite[Section 13]{Bi}):
$$
\bigl(\rho^{(n)}_t, g_1\bigr)\bigl(\rho^{(n)}_t, g_2\bigr)\cdot\ldots\cdot\bigl(\rho^{(n)}_t, g_m\bigr);\quad n = 1, 2, \ldots
$$
Therefore, for every $\Phi \in \mathfrak A$, the following sequence is tight in $D[0, T]$:
$$
\Phi\bigl(\rho^{(n)}_t\bigr),\, t \in [0, T];\quad n = 1, 2, \ldots
$$
Apply criteria of relative compactness: \cite[Proposition 3.9.1]{EthierKurtz}, and complete the proof.

\subsection{Characterization of weak limits} In this step, we will characterize any weak limit point $\rho = (\rho_t,\, 0 \le t \le T)$ of $(\rho^{(n)})$. We shall think in terms of cumulative distribution functions: Let $(\rho^{(n_k)})$ be any subsequence weakly converging in $D([0, T], \mathcal P_q)$ to $\rho$, and let $F_k(t, \cdot)$ be the cumulative distribution function of $\rho^{(n_k)}_t$. Without loss of generality, by the Skorohod representation theorem we can assume convergence a.s. in $D([0, T], \mathcal P_q)$.

The standard approach to derive the limit is to adapt the arguments in \cite[Lemma 1.5]{Jo}; however, the arguments in \cite{JR} cannot be extended to prove uniqueness of solutions for \eqref{SPDE2}. We adopt a different and a much simpler approach to derive the limit and to prove uniqueness of limits. We again use the idea from subsection 2.1:  Reduce the particle system in \eqref{main_sde} to the rank-based system, with $\gamma = 0$. We use the notation from there. The arrow $\stackrel{\mathbb{P}}{\longrightarrow}$ indicates convergence in probability. 

Under Assumption~\ref{main_asmp}, the Cauchy problem~\eqref{PME} admits a unique solution $R$ (in the distributional sense) with distributional derivative $R_x$, which is a classic function, and
\begin{equation}
\label{eq:bdd-Rx}
C_* := \sup\limits_{t \in [0, T]}\sup\limits_{x \in \mathbb R}|R_x(t, x)|.
\end{equation}
 We claim that:
\begin{equation}\label{uniform_convergence}
\sup\limits_{t \in [0,T]}\sup\limits_{x \in \rr}\Big|F_k(t, x)-R\Big(t,x-\int_0^t \gamma(\rho(s))\mathrm{d}W(s)\Big)\Big| \stackrel{\mathbb{P}}{\longrightarrow}0,\quad k \to \infty.
\end{equation} 
Similarly to~\eqref{main_eqn}, we use the following notation for shorthand:
$$
\Gamma(t) := \int_0^t \gamma(\rho(s))\,\mathrm{d}W(s),\quad \Gamma_k(t) := \int_0^t \gamma(\rho^{(n_k)}(s))\mathrm{d}W(s)
$$
In view of \eqref{eq:Y} and {\eqref{observations}, we need only to show convergence in probability:
\begin{equation}\label{uni_conv}
\sup\limits_{t \in [0,T]}\sup\limits_{x \in \rr}\big|F_{\mu^{(n_k)}(t)}(x-\Gamma_k(t))-R(t,x-\Gamma(t))\big|\stackrel{\mathbb{P}}{\longrightarrow}0. 
\end{equation}
We apply the triangle inequality to bound from above the left-hand side of~\eqref{uni_conv}:
\begin{align}
\label{eq:triangle-ineq}
\begin{split}
&\sup\limits_{t \in [0,T]}\sup\limits_{x \in \rr}\big|F_{\mu^{(n_k)}(t)}(x-\Gamma_k(t))- R(t,x-\Gamma_k(t))\big|
\\ & +\sup\limits_{t \in [0,T]}\sup\limits_{x \in \rr}\big|R(t,x-\Gamma_k(t)) - R(t, x - \Gamma(t))\big| \\ 
& \le \sup\limits_{t \in [0,T]}\sup\limits_{x\in \rr}\big|F_{\mu^{(n_k)}(t)}(x) - R(t,x)\big|
+ C_*\sup\limits_{t \in [0,T]}\big|\Gamma_k(t) - \Gamma(t)\big|.
\end{split}
\end{align}
The claim and its proof in proposition \cite[Equation 5.17, Proposition 5.1]{KoS} imply:
$$
\sup\limits_{t \in [0,T]}\sup\limits_{x\in \rr}\big|F_{\mu^{(n_k)}(t)}(x) - R(t,x)\big|\stackrel{\mathbb{P}}{\longrightarrow}0.
$$
It remains to show that 
\begin{equation}
\label{eq:conv-in-prob}
\sup_{t \in [0,T]}|\Gamma_k(t) - \Gamma(t)| \stackrel{\mathbb{P}}{\longrightarrow} 0.
\end{equation}
In light of Assumption~\ref{main_asmp} (c), we can estimate
\begin{equation}
\label{eq:inter}
\langle\Gamma_k-\Gamma\rangle_T = \int_0^T\left(\gamma(\rho_{n_k})(t)) - \gamma(\rho(t))\right)^2\,\mathrm{d}t \le L_{\gamma}^2\int_0^T\mathcal W^2_1(\rho^{(n_k)}(t), \rho(t))\,\mathrm{d}t.
\end{equation}
We have the following convergence in law, and therefore in probability:
\begin{equation}
\label{eq:need}
\int_0^T\mathcal W^2_1\left(\rho^{n_k}(s), \rho(s)\right)\,\mathrm{d}s \le T\cdot\sup\limits_{0 \le t \le T}\mathcal W_1^2\left(\rho^{n_k}(s), \rho(s)\right) \to 0
\end{equation}
Combining~\eqref{eq:inter} and~\eqref{eq:need} with \cite[Chapter 1, Problem 5.25]{KS}, we prove~\eqref{eq:conv-in-prob}. This completes the proof of~\eqref{uni_conv}. Next, 
\begin{equation}
\label{eq:int-conv}
\sup\limits_{0 \le t \le T}\int_{\rr}|F_k(t, x) - F_{\rho(t)}(x)| \,\mathrm{d}x = \sup\limits_{0 \le t \le T}\mathcal W_1(\rho^{(n_k)}(t), \rho(t)) \to 0
\end{equation}
in law (and therefore in probability) as $k \to \infty$. Combining~\eqref{uni_conv} with~\eqref{eq:int-conv}, we get: for any bounded interval $I \subseteq \mathbb R$, almost surely,
$$
\sup\limits_{0 \le t \le T}\int_I\Big|F_{\rho(t)}(x) - R\Big(t,x-\int_0^t \gamma(\rho(s))\mathrm{d}W(s)\Big)\Big|\,\mathrm{d}x = 0.
$$
Thus almost surely for every $(t, x) \in [0, T]\times I$ we get~\eqref{main_eqn}. Representing the real line as a countable union of such intervals $I$, and noting that intersection of countably many almost sure events is also almost sure, we prove~\eqref{main_eqn} a.s. for all $(t, x) \in [0, T]\times\mathbb R$. 

\subsection{Derivation of the SPDE} Assumption 1.1(e) with \cite[Lemma 2.7]{JR} yield classical regularity for $R$. Apply It\^o's formula to $G(t, x) = F_{\rho(t)}(x)$ in \eqref{main_eqn}: 
\begin{align*}
\mathrm{d} G & = \frac{\partial R}{\partial t}(t,x-\Gamma(t))\,\mathrm{d}t - \frac{\partial R}{\partial x}(t,x-\Gamma(t))\,\gamma(\rho(t))\,\mathrm{d}W(t)  + \frac12\frac{\partial^2R}{\partial x^2}(t,x-\Gamma(t))\,\gamma^2(\rho(t))\,\mathrm{d}t \\&=  -B(R(t,x-\Gamma(t)))_x\,\mathrm{d}t+\Sigma(R(t,x-\Gamma(t)))_{xx}\,\mathrm{d}t  - \frac{\partial R}{\partial x}(t,x-\Gamma(t))\,\gamma(\rho(t))\,\mathrm{d}W(t) \\ & \quad + \frac12\frac{\partial^2R}{\partial x^2}(t,x-\Gamma(t))\,\gamma^2(\rho(t))\,\mathrm{d}t,
\end{align*}
where the last equality is a consequence of \eqref{PME}. Noting that 
$$
G_x(t,x) = R_x(t,x-\Gamma(t))\quad \mbox{and}\quad G_{xx}(t,x) = R_{xx}(t,x-\Gamma(t)),
$$
we obtain the SPDE \eqref{SPDE2}. In the fourth and final section, we prove  that the solution to this functional equation~\eqref{main_eqn} is unique. Taken together, all of this proves Theorem~\ref{thm:lln} with convergence in the Skorohod space instead of the uniform convergence. Since the corresponding measure-valued process $\rho$ is a.s. continuous with respect to time $t$, it is an element of $C([0, T], \mathcal P_q)$. The same can be said about the pre-limit processes $\rho^{(n)}$, and thus the convergence takes place in $C([0, T], \mathcal P_q)$: See \cite[Chapter 12]{Bi}.

\subsection{Uniqueness of  the limit} Let $\rho_1$ and $\rho_2$ be in $C([0, T], \mathcal P_p(\rr))$ with continuous cumulative distribution functions $F_i(t, x) := F_{\rho_i(t)}(x)$, $i = 1, 2$, satisfying~\eqref{main_eqn}. Denote $\Gamma_i(t) := \int_0^t\gamma(\rho_i(s))\,\mathrm{d}W(s)$ for $i = 1, 2$, we get:
\begin{align}
\label{eq:diff-F}
\begin{split}
F_1(t, x) - F_2(t, x) & = R(t,x-\Gamma_1(t)) - R(t,x-\Gamma_2(t)) \\ & = \int_0^1R_x\left[t, x - \Gamma_1(t)\theta - \Gamma_2(t)(1 - \theta)\right]\,\mathrm{d}\theta\cdot(\Gamma_1(t) - \Gamma_2(t)).
\end{split}
\end{align}
We can represent the Wasserstein distance between $\rho_1$ and $\rho_2$ as follows \cite[p.64]{SW}:
\begin{equation}
\label{eq:W-1}
\mathcal W_1(\rho_1(t), \rho_2(t))  = \int_\rr|F_1(t, x) - F_2(t, x)|\,\mathrm{d}x.
\end{equation}
Applying~\eqref{eq:bdd-Rx} above to~\eqref{eq:diff-F} with~\eqref{eq:W-1} and interchanging integrations by Fubini's theorem, we obtain
\begin{align}
\label{eq:W-1-est}
\mathcal W_1(\rho_1(t), \rho_2(t))  \le \int_\rr R_x(t,x)\,\mathrm{d}x \cdot\Bigl|\int_0^t\gamma(\rho_1(s))-\gamma(\rho_2(s))\,\mathrm{d}W(s)\Bigr|.
\end{align}
Note that $R_x(t, \cdot)$ is the probability density function, which integrates to 1. We square both sides in~\eqref{eq:W-1-est}, take expectation, and apply the Doob's martingale inequality:
\begin{equation*}\label{eqn}
\begin{split}
\E&\Big[\mathcal W_1^2(\rho_1(t),\rho_2(t))\Big] \le \,4\E\int_0^t\big|\gamma(\rho_1(s))-\gamma(\rho_2(s))\big|^2\,\mathrm{d}s \\ &\le 4L^2_{\gamma}\,\int_0^t \E\,\mathcal W_1^2(\rho_1(s),\rho_2(s))\,\mathrm{d}s \le 4L_{\gamma}^2\int_0^t\E\,\mathcal W_1^2(\rho_1(s), \rho_2(s))\,\mathrm{d}s.
\end{split}
\end{equation*}
Gronwall's lemma implies uniqueness. 





\begin{thebibliography}{BFPS}

\bibitem[BFK]{BFK} A.~D.~Banner, E.~R.~Fernholz, I.~Karatzas (2005). Atlas models of equity markets. \textit{Ann. Appl. Probab.} \textbf{15}, 2296--2330.


\bibitem[BP]{BP} R.~F.~Bass, E.~Pardoux (1987). Uniqueness for diffusions with piecewise constant coefficients. \textit{Probab. Theory Related Fields} \textbf{76} 557--572. 

\bibitem[Bi]{Bi} P.~Billingsley (1999). \textit{Convergence of probability measures}. 2nd edition. Wiley.

\bibitem[CP]{CP} S.~Chatterjee, S.~Pal (2011). A combinatorial analysis of interacting diffusions. \textit{J. Theoret. Probab.} \textbf{24}, pp. 939--968.

\bibitem[DSVZ]{DSVZ} A.~Dembo, M.~Shkolnikov, S.~R.~S.~Varadhan, O.~Zeitouni (2016). Large deviations for diffusions interacting through their ranks. \textit{Comm. Pure Appl. Math.} \textbf{69}, 1259--1313.


\bibitem[EK]{EthierKurtz} S.~N.~Ethier, T.~G.~Kurtz (2005). \textit{Markov Processes: Characterization and Convergence.} Wiley.

\bibitem[Fa]{Fan} J.~Fan, Y.~Liao, H.~Liu (2016). An overview on the estimation of large covariance and precision matrices. \textit{Economet. J.} \textbf{16} C1--C32. 

\bibitem[Fo]{Folland} G.~B.~Folland (1999). \textit{Real Analysis: Modern Techniques and Their Applications.} 2nd edition, Wiley.

\bibitem[Fe]{Fe} E.~R.~Fernholz (2002). \textit{Stochastic portfolio theory}. Applications of Mathematics \textbf{48}. Springer.

\bibitem[FK]{FK} R.~Fernholz, I.~Karatzas (2009). Stochastic portfolio theory: an overview. In: A.~Bensoussan, Q.~Zhang (eds.) \textit{Handbook of Numerical Analysis}. Mathematical
Modeling and Numerical Methods in Finance \textbf{XV}, pp. 89--167. North-Holland, Oxford.

\bibitem[Ga]{Ga} J.~G\"artner (1988). On the McKean-Vlasov limit for interacting diffusions. \textit{Math. Nachr.} \textbf{137}, 197--248.

\bibitem[Gi]{Gi} B.~H.~Gilding (1989). Improved theory for a nonlinear degenerate parabolic equation. \textit{Ann. Sc. Norm. Super. Pisa Cl. Sci.} \textbf{16}, 165--224.

\bibitem[GS]{GS} B.~Gess, P.~Souganidis (2017). Stochastic non-isotropic degenerate parabolic-hyperbolic equations. \textit{Stoch. Proc. Appl.} \textbf{127}, 2961--3004.

\bibitem[IPS]{IPS} T.~Ichiba, S.~Pal, M.~Shkolnikov (2013). Convergence rates for rank-based models with applications to portfolio theory. \textit{Probab. Theory Related Fields} \textbf{156}, 415--448.

\bibitem[Jo]{Jo} B.~Jourdain (2000). Diffusion processes associated with nonlinear evolution equations for signed measures. \textit{Methodol. Comput. Appl. Probab.} \textbf{2}, 69--91. 

\bibitem[JR]{JR} B.~Jourdain, J. Reygner (2013). Propagation of chaos for rank-based interacting diffusions and long time behaviour of a scalar quasilinear parabolic equation. \textit{Stochastic Partial Differential Equations: Analysis and Computations} \textbf{1}, 455--506.

\bibitem[KS]{KS} I.~Karatzas, S.~Shreve (1991). \textit{Brownian motion and stochastic calculus}. 2nd edition. Springer.

\bibitem[KoS]{KoS} P.~Kolli, M.~Shkolnikov (2018). SPDE limit of the global fluctuations in rank-based models. \textit{Ann. Probab.} \textbf{46}, 1042-1069. 

\bibitem[LSU]{LSU} O.~A.~Ladyzenskaja, V.~A.~Solonnikov, N.~Uralceva (1968). \textit{Linear and quasilinear equations of parabolic type}. Translations of Mathematical Monographs \textbf{23}. American Mathematical Society.

\bibitem[ILS]{MyOwn19} T.~Ichiba, M.~Ludkovski, A.~Sarantsev (2018). Dynamic contagion in a banking system with births and defaults. Available at arXiv:1807.09897.

\bibitem[LPS]{LPS}P.~L.~Lions, B.~Perthame, P.~Souganidis (2013). Scalar conservation laws with rough (stochastic) fluxes.  \textit{Stochastic Partial Differential Equations: Analysis and Computations} \textbf{1}, 664-686

\bibitem[PO]{Cov} M.~Pourahmadi (2013). \textit{High-Dimensional Covariance Estimation.} Wiley. 

\bibitem[S]{S} M.~Shkolnikov (2012). Large systems of diffusions interacting through their ranks. \textit{Stoch. Proc. Appl.} \textbf{122}, 1730--1747.

\bibitem[SW]{SW} G.~R.~Shorack, J.~A.~Wellner (1986). Empirical processes with applications to statistics. Wiley.

\bibitem[SV]{SV} D.~W.~Stroock, S.~R.~S.~Varadhan (2006). \textit{Multidimensional diffusion processes}. Springer.




\end{thebibliography}
\end{document}